%% file: 0oricplxcob.tex
\documentclass[11pt]{amsart}
\usepackage{amssymb,amscd}
\usepackage{srcltx}
\usepackage{graphicx}
\usepackage[colorlinks=true, urlcolor=blue,bookmarks=true,bookmarksopen=true, citecolor=blue,hypertex]{hyperref}

\numberwithin{equation}{section}

\newtheorem{defn}{Definition}[section]
\newtheorem{theorem}{Theorem}[section]

\newtheorem{example}[theorem]{Example}

\newtheorem{lemma}[theorem]{Lemma}
\newtheorem{prop}[theorem]{Proposition}

\def \begineq{\begin{equation}}
\def \endeq{\end{equation}}

\def \bb{\mathbb}

\def \CC{{\bb{C}}}
\def \CP{{\bb{CP}}}

\def \RR{{\bb{R}}}
\def \TT{{\bb{T}}}

\def \ZZ{{\bb{Z}}}

\def \({\left(}
\def \){\right)}
\def \<{\langle}
\def \>{\rangle}

\begin{document}

\title{ Oriented Cobordism of $\CP^{2k+1}$}

\author[S. Sarkar]{Soumen Sarkar}
\address{\it{Theoretical Statistics and Mathematics Unit, Indian
Statistical Institute, 203 B. T. Road, Kolkata 700108, India}}

\email{soumens$_-$r@isical.ac.in}



\abstract We give a new construction of oriented manifolds having the boundary $\CC P^{2k+1}$
for each $k \geq 0$. The main tool is the theory of quasitoric manifolds.
\endabstract
\maketitle

\section{Introduction}\label{intro}
Cobordism was explicitly introduced by Lev Pontryagin in geometric work on manifolds.
In the early 1950's Ren\'{e} Thom showed that cobordism groups could be computed by
results of homotopy theory. In \cite{Tho} Thom showed that the cobordism classes
of $G$-manifolds, for a Lie group $G$, are in one to one correspondence with the
elements of the homotopy group of the Thom space of the group $G \subseteq O_m$.
In this article, for each $k \geq 0$, we construct oriented manifold
with boundary whose boundary is $\CP^{2k+1}$. To construct these manifold
with boundary we use the theory of quasitoric manifolds.

Quasitoric manifolds were introduced by Davis and Januskiewicz in \cite{DJ}.
A manifold with quasitoric boundary is a manifold with boundary where the boundary is a disjoint union of some quasitoric
manifolds. The strategy of our proof is to first construct some compact orientable manifolds with quasitoric boundary.
Then identifying suitable boundary components using certain equivariant homeomorphisms
we obtain oriented manifolds with the boundary $\CP^{2k+1}$ for each $k \geq 0$, see theorem \ref{gr}.

\section{Quasitoric manifolds}
A $k$-dimensional simple polytope is a convex polytope where exactly $k$ bounding hyperplanes meet at 
each vertex. The codimension one faces of a convex polytope are called facets.
Let $\mathcal{F}(P)$ be the set of facets of a $k$-dimensional simple polytope $P$.
Following \cite{BP} we give definition of quasitoric manifold, characteristic function and classification.
\begin{defn}
An action of $\TT^k$ on a $2k$-dimensional manifold $M$ is said to be locally
standard if every point $y \in M $ has a $\TT^k$-stable open neighborhood
$U_y$ and a homeomorphism $\psi : U_y \to V$, where $V$ is a $\TT^k$-stable
open subset of $\CC^k$, and an isomorphism $\delta_y : \TT^k \to \TT^k$ such that
\ $\psi (t\cdot x) = \delta_y (t) \cdot \psi(x)$ for all $(t,x) \in \TT^k \times U_y$.
\end{defn}

\begin{defn}\label{qtd02}
A $\TT^k$-manifold $M$ is called a quasitoric manifold over $P$ if the following conditions are satisfied: 
\begin{enumerate}
\item the $\TT^k$ action is locally standard,
\item there is a projection map $\mathfrak{q}: M \to P$ constant on $\TT^k$ orbits which maps every
$l$-dimensional orbit to a point in the interior of a codimension-$l$ face of $P$.
\end{enumerate}
\end{defn}

All complex projective spaces $\CP^{k}$ and their equivariant connected sums, products are quasitoric manifolds.
\begin{lemma}[\cite{DJ}, Lemma 1.4]\label{clema}
Let $\mathfrak{q}: M \to P$ be a $2k$-dimensional quasitoric manifold over $P$. There is a projection map
$f: \TT^k \times P \to M$ so that for each $q \in P$, $f$ maps $\TT^k \times q$ onto $ \mathfrak{q}^{-1}(q)$.
\end{lemma}

A quasitoric manifold $M$ over $P$ is simply connected. So $M$ is orientable.
A choice of orientation on $\TT^k $ and $ P$ gives an orientation on $M$.

Define an equivalence relation $\sim_2$ on $\ZZ^k$ by $x \sim_2 y$ if and only if $y = -x$.
Denote the equivalence class of $x$ in the quotient space $\ZZ^k/\ZZ_2$ by $[x]$.
\begin{defn} 
A function $\lambda : \mathcal{F}(P) \to \ZZ^k/\ZZ_2$ is called characteristic function if the submodule generated by
$\{\lambda(F_{j_1}), \ldots, \lambda(F_{j_l})\}$ is an $l$-dimensional direct summand of $\ZZ^k$ whenever
the intersection of the facets $F_{j_1}, \ldots, F_{j_l}$ is nonempty.

The vectors $ \lambda(F_{j})$ are called characteristic vectors and the pair $(P, \lambda)$ is called a characteristic pair.
\end{defn}

In \cite{DJ} the authors show that we can construct a quasitoric manifold from the pair $(P, \lambda)$.
Also given quasitoric manifold we can associate a characteristic pair to it up to choice of signs of
characteristic vectors.

\begin{defn} 
Let $\delta : \TT^k \to \TT^k$ be an automorphism. Two quasitoric manifolds $M_1$ and $M_2$ over the same
polytope $P$ are called $\delta$-equivariantly homeomorphic if there is a homeomorphism $ f : M_1 \to M_2$
such that $f(t \cdot x) = \delta(t)\cdot f(x)$ for all $(t, x ) \in \TT^k \times M_1$.
\end{defn}
The automorphism $\delta$ induces an automorphism $\delta_{\ast}$ of the poset of subtori of $\TT^k$.
This automorphism descends to a $\delta$-$translation$ of characteristic pairs, in which the two
characteristic functions differ by $\delta_{\ast}$. Using Lemma \ref{clema} we can prove the following Proposition.

\begin{prop}[\cite{BP}, Proposition 5.14]\label{probi}
There is a bijection between $\delta$-equivariant homeomorphism classes of quasitoric manifolds and
$\delta$-translations of characteristic pairs $(P, \lambda)$.
\end{prop}
Suppose $\delta$ is the identity automorphism of $\TT^k$. From Proposition \ref{probi} we have two quasitoric
manifolds are equivariantly homeomorphic if and only if their characteristic functions are the same.

\section{Some manifolds with quasitoric boundary}\label{tcobdd}
Set $n=2(k+1)$. Corresponding to each even $k \geq 0$ we construct a manifold with quasitoric boundary.
Let $ \{ A_j: j=0, \ldots, n \} $ be the standard basis of $\RR^{n+1}$. Let
\begin{equation}\label{eq01}
 \bigtriangleup^n = \{(x_0, x_1, \ldots, x_n) \in \RR^{n+1} : ~ x_j \geq 0 ~\mbox{and} ~ \Sigma_{0}^{n}x_j = 1\}.
\end{equation}
Then $ \bigtriangleup^n $ is an $n$-dimensional simplex with vertices $ \{ A_j: j=0, \ldots, n \} $ in $\RR^{n+1}$. Define
\begin{equation}
 \bigtriangleup_j^{n-1} = \{(x_0, x_1, \ldots, x_n) \in \bigtriangleup^n : x_j =0\}.
\end{equation}
So $ \bigtriangleup_j^{n-1} $ is the facet of $ \bigtriangleup^n$ not containing the vertex $A_j$.
Let $F$ be the largest face of $\bigtriangleup^n$ not containing the vertices $A_{j_1}, \ldots, A_{j_l}$. Then
\begin{equation}
 F= \cap_{i=1}^{l}\bigtriangleup_{j_i}^{n-1} = \{(x_0, x_1, \ldots, x_n) \in \bigtriangleup^n : x_{j_i} =0,~ i=1, \ldots, l\}.
\end{equation}

 Define a function $ \eta :\{ \bigtriangleup_j^{n-1} : j=0, \ldots, n \} \longrightarrow \ZZ^{n-1} $
as follows.

\begin{equation}
\eta(\bigtriangleup_{n-j}^{n-1}) = \left\{ \begin{array}{ll} (0, \ldots, 0, 1,0, \ldots, 0) & \mbox{if}~ 0 \leq j < \frac{n}{2}-1, ~
\mbox{here} ~ 1 ~ \mbox{is in the} ~ (j+1)\mbox{-th place} \\
 (1, \ldots, 1, 1, 0, \ldots, 0) & \mbox{if} ~ j = \frac{n}{2}-1, ~ \mbox{here} ~ 1 ~ \mbox{occurs up to} ~ \frac{n}{2}\mbox{-th place}\\
 (0, \ldots, 0, 1,0, \ldots, 0) & \mbox{if} ~ \frac{n}{2} \leq j < n, ~ \mbox{here} ~ 1~ \mbox{is in the} ~ j\mbox{-th place}\\
 (0, \ldots, 0, 1, 1, \ldots, 1) & \mbox{if}~ j = n, ~ \mbox{here} ~ 0 ~ \mbox{occurs up to} ~ (\frac{n}{2}-1)\mbox{-th place}.
\end{array} \right.
\end{equation}

Define\label{equxi}
\begin{equation}
 \eta_j := \eta(\bigtriangleup_{n-j}^{n-1}), ~ \mbox{for all} ~ j = 0, 1, \ldots, n.
\end{equation}

\begin{example}\label{exam1}
For $n=4$, let $ \bigtriangleup^4 $ be the $4$-simplex in $\RR^5$ with vertices $A_0, A_1, A_2, A_3,$ $ A_4$ (see figure \ref{egch060}).
Define a function $\eta$ from the set of facets of $ \bigtriangleup^4$ to $\ZZ^3$ by,
\begin{equation}
\eta(\bigtriangleup_{4-j}^3)= \left\{ \begin{array}{ll} (1,0,0) & \mbox{if}~  j= 0\\
 (1,1,0) & \mbox{if}~ j = 1\\
 (0,1,0) & \mbox{if} ~ j= 2\\
 (0,0,1) & \mbox{if} ~ j= 3\\
 (0,1,1) & \mbox{if}~ j = 4.
 \end{array} \right.
\end{equation}
\end{example}

\begin{figure}[ht]
\centerline{
 \scalebox{0.70}{
  \input{egch060.pstex_t}
          }
 }
\caption {The $4$-simplex  $ \bigtriangleup^4$.}
\label{egch060}
\end{figure}

Suppose the faces $F^{\prime}$ and $F^{\prime \prime}$ of $\bigtriangleup^n$ are the intersection of facets
$\{ \bigtriangleup_{n}^{n-1}, \bigtriangleup_{n-1}^{n-1}, \ldots,$ $\bigtriangleup_{\frac{n}{2}}^{n-1}\}$ and
$ \{ \bigtriangleup_{\frac{n}{2}}^{n-1}, \bigtriangleup_{\frac{n}{2}-1}^{n-1}, \ldots, \bigtriangleup_{0}^{n-1} \}$
respectively. Then
\begin{equation}
 F^{\prime}= \{(x_0, x_1, \ldots, x_n) \in \bigtriangleup^n :  x_{\frac{n}{2}} =0, \ldots, x_n =0 \},
\end{equation}
\begin{equation}
 F^{\prime\prime}= \{(x_0, x_1, \ldots, x_n) \in \bigtriangleup^n : x_0 =0, \ldots, x_{\frac{n}{2}} = 0\}.
\end{equation}
Hence dim$(F^{\prime})$ $=$ dim$(F^{\prime \prime}) = \frac{n}{2}-1 \geq 1$.
The set of vectors $\{ \eta_0, \ldots, \eta_{\frac{n}{2}}\}$ and $\{\eta_{\frac{n}{2}},  \ldots, \eta_n\}$ 
are linearly dependent sets in $\ZZ^{n-1}$. But the submodules generated by the vectors
$\{ \eta_0, \ldots, $ $\widehat{\eta_{j}}, \ldots, \eta_{\frac{n}{2}}\}$ and
$\{\eta_{\frac{n}{2}}, \eta_{\frac{n}{2}+1}, \ldots, \widehat{\eta_{l}}, \ldots, \eta_n\}$ are $\frac{n}{2}$-dimensional
direct summands of $\ZZ^{n-1}$ for each $j= 0, \ldots, \frac{n}{2}$ and
$l=\frac{n}{2}, \ldots, n$ respectively. Here the symbol $~\widehat{}~$ represents the omission of the corresponding entry.

Suppose $e$ is an edge of $\bigtriangleup^n$ not contained in $F^{\prime} \cup F^{\prime \prime}$. Then 
$e= \cap_{_{j=1}}^{^{n-1}} \bigtriangleup_{n-{l_j}}^{n-1}$ for some $\{l_j : j=1, \ldots, n-1\} \subset \{0, 1, \ldots, n\}$.
Observe that $\{ \eta_0, \ldots, \eta_{\frac{n}{2}}\} \nsubseteq \{ \eta_{l_1}, \ldots, \eta_{l_{n-1}}\}$ and
$\{\eta_{\frac{n}{2}}, \eta_{(\frac{n+2}{2})}, \ldots, \eta_n\} \nsubseteq \{ \eta_{l_1}, \ldots, \eta_{l_{n-1}}\}$.
Hence the set of vectors $\{ \eta_{l_1}, \ldots, \eta_{l_{n-1}}\}$ form a basis of $\ZZ^{n-1}$.

Let $r_1$ be a positive real number such that $ r_1 <\frac{1}{4}$.
Consider the following affine hyperplanes in $\RR^{n+1}$.
\begin{equation}
H_1= \{(x_0, x_1, \ldots, x_n) \in \RR^{n+1} :~ x_{\frac{n}{2}} + \cdots + x_n =r_1 \}.
\end{equation}
\begin{equation}
 H_2= \{(x_0, x_1, \ldots, x_n) \in \RR^{n+1} :~ x_0 + \cdots + x_{\frac{n}{2}} = r_1\}.
\end{equation}
\begin{equation}
 H_3= \{(x_0, x_1, \ldots, x_n) \in \RR^{n+1} : ~x_{\frac{n}{2}} = 1-r_1\}.
\end{equation}
\begin{equation}
 H = \{(x_0, x_1, \ldots, x_n) \in \RR^{n+1} : ~x_0 + \cdots + x_n = 1\}.
\end{equation}
\begin{figure}[ht]
\centerline{
\scalebox{0.90}{
\input{egch06.pstex_t}
 }
 }
\caption {The $4$-simplex $ \bigtriangleup^4$ and the affine hyperplanes $H_1, H_2 $ and $ H_3$.}
\label{egch061}
\end{figure}
Then $ \bigtriangleup^n \subset H$ and the intersections
$ \bigtriangleup^n \cap H_1 \cap H_2, ~ \bigtriangleup^n \cap H_1 \cap H_3, ~ \bigtriangleup^n \cap H_3 \cap H_2$ are empty.
We cut off an open neighborhood of faces $F^{\prime}$, $F^{\prime \prime}$ and $\{A_{\frac{n}{2}}\}$
by affine hyperplanes $H_1\cap H$, $H_2 \cap H$ and $ H_3 \cap H$ respectively in $H$.
Let $H_j^{\prime}$ be the closed half space associated to the affine hyperplane $ H_j$ such that
the interior of half spaces $H_1^{\prime}, H_2^{\prime}, H_3^{\prime}$ do not contain the
faces $F^{\prime}, F^{\prime \prime}, \{A_{\frac{n}{2}}\}$ respectively.
We illustrate such hyperplanes for the case $n=4$ in figure \ref{egch061}. Define
\begin{equation}
 {\bigtriangleup^n_Q} = \bigtriangleup^n \cap H_1^{\prime} \cap H_2^{\prime} \cap H_3^{\prime}, ~ P_1 = \bigtriangleup^n \cap H_1, ~
P_2 = \bigtriangleup^n \cap H_2 ~ \mbox{and} ~ P_3 = \bigtriangleup^n \cap H_3.
\end{equation}

The convex polytope ${\bigtriangleup^n_Q}$ is a simple convex polytope of dimension $n$
and the polytopes $P_1$, $P_2$ and $P_3$ are also facets of ${\bigtriangleup^n_Q}$.
The polytopes $P_1$, $P_2$ and $P_3$ are given by the following equations. 
\begin{equation}\label{equp1}
P_1= \{(x_0, x_1, \ldots, x_n) \in \bigtriangleup^n :~ x_0 + \cdots + x_{\frac{n}{2}-1} = 1-r_1 ~ \mbox{and}
~ x_{\frac{n}{2}} + \cdots + x_n =r_1 \}.
\end{equation}
\begin{equation}\label{equp2}
 P_2= \{(x_0, x_1, \ldots, x_n) \in \bigtriangleup^n :~ x_0 + \cdots + x_{\frac{n}{2}} = r_1  ~ \mbox{and}
~ x_{\frac{n}{2}+1} + \cdots + x_n = 1-r_1 \}.
\end{equation}
\begin{equation}\label{equp3}
 P_3= \{(x_0, x_1, \ldots, x_n) \in \bigtriangleup^n : ~x_{\frac{n}{2}} = 1-r_1 ~ \mbox{and}
~ x_0 + \cdots + \widehat{x}_{\frac{n}{2}} + \cdots + x_n =r_1 \}.
\end{equation}

By equations \ref{equp1} and \ref{equp2}, the convex polytopes $P_1$ and $P_2$ are diffeomorphic
to the product $\bigtriangleup^{\frac{n}{2}-1} \times \bigtriangleup^{\frac{n}{2}}$. From equation
\ref{equp3} $P_3$ is diffeomorphic to the simplex $\bigtriangleup^{n-1}$. The facets of $P_1$, $P_2$
and $P_3$ are given by the following equations respectively.
\begin{equation}\label{equfp1}
\bigtriangleup_j^{n-1} \cap P_1= \{(x_0, x_1, \ldots, x_n) \in P_1 : x_j =0\} ~ \mbox{for all} ~ j\in \{0, \ldots, n\}.
\end{equation}
\begin{equation}\label{equfp2}
 \bigtriangleup_j^{n-1} \cap P_2= \{(x_0, x_1, \ldots, x_n) \in P_2 : x_j = 0\} ~ \mbox{for all} ~ j\in \{0, \ldots, n\}.
\end{equation}
\begin{equation}
 \bigtriangleup_j^{n-1} \cap P_3= \{(x_0, x_1, \ldots, x_n) \in P_3 : x_j =0\}~ \mbox{for all}~ j \in \{0, \ldots,
\widehat{\frac{n}{2}}, \ldots, n\}.
\end{equation}

Now we want to construct $(2n-1)$-dimensional manifold with quasitoric boundary.
Let $F$ be a face of $\bigtriangleup^n$ of codimension $l$.
Then $F = \bigtriangleup^{n-1}_{j_1} \cap \ldots \cap \bigtriangleup^{n-1}_{j_l}$ for a unique
$\{j_1, \ldots, j_l\} \subset \{0,1, \ldots, n\}$. Suppose $\TT_F$ be the torus subgroup of
$\TT^{n-1}$ determined by the submodule generated by $\{ \eta_{j_1}, \ldots, \eta_{j_l}\}$ in $\ZZ^{n-1}$.
Assume $\TT_{\bigtriangleup^n} = \{1\}$.
We define an equivalence relation $\sim_{\eta}$ on the product $\TT^{n-1} \times \bigtriangleup^n$ as follows,
\begin{equation}\label{defeqiv}
(s,p) \sim_{\eta} (t,q) ~ \mbox{if and only if}~  p = q ~ \mbox{and} ~ ts^{-1} \in \TT_{F}
\end{equation}
where $ F \subset \bigtriangleup^n $ is the unique face containing the point $ p $ in its relative interior.
Restrict the equivalence relation $\sim_{\eta}$ on $\TT^{n-1} \times {\bigtriangleup^n_Q}$.
Define $$W(\bigtriangleup^n_Q, \eta) := (\TT^{n-1} \times {\bigtriangleup^n_Q})/\sim_{\eta}$$ to be the quotient space.
So $W(\bigtriangleup^n_Q, \eta)$ is a $\TT^{n-1}$-space. The fixed point set of this action corresponds bijectively
to $\bigtriangleup^n_Q \cap \bigtriangleup^n(1)$, where $\bigtriangleup^n(1)$ is the $1$-skeleton of $\bigtriangleup^n$.
Let $$\mathfrak{p} : W(\bigtriangleup^n_Q, \eta) \to \bigtriangleup^n_Q,$$ defined by
$\mathfrak{p}([s,p]^{\sim_{\eta}}) = p$, be the corresponding orbit map.

Let $\eta^{1}$, $\eta^{2}$ and $\eta^{3}$ be the restriction of the function $\eta$
on the set of facets of $P_1$, $P_2$ and $P_3$ respectively. Define
\begin{equation}\label{char1}
\eta^i_j := \eta^i(\bigtriangleup_{n-j}^{n-1} \cap P_i)= \left\{ \begin{array}{ll} \eta_j & \mbox{if}~  i= 1, 2 ~\mbox{and}~ 
j\in \{0,1, \ldots, n\} \\
 \eta_j & \mbox{if}~ i = 3 ~\mbox{and}~ j\in \{0,1, \ldots, \widehat{\frac{n}{2}}, \ldots, n\}.
 \end{array} \right.
\end{equation}

Let $v$ be a vertex of $P_i$. So $v$ belongs to the relative interior of a unique edge $e_v$ of $\bigtriangleup^n$
not contained in $F^{\prime} \cup F^{\prime \prime}$. If $e_v = \cap_{_{j=1}}^{^{n-1}} \bigtriangleup_{n-{l_j}}^{n-1}$ for
some $\{l_j : j=1, \ldots, n-1\} \subset \{0, 1, \ldots, n\}$, the vectors $\{ \eta_{l_1}, \ldots, \eta_{l_{n-1}}\}$
form a basis of $\ZZ^{n-1}$. So $$v = \cap_{_{j=1}}^{^{n-1}} (\bigtriangleup_{l_j}^{n-1} \cap P_i)$$ and
the vectors $\{ \eta^i_{l_1}, \ldots, \eta^i_{l_{n-1}}\}$ form a basis of $\ZZ^{n-1}$.
So $\eta^{i}$ defines the characteristic function of a quasitoric manifold $M(P_i, \eta^i)$ over $P_i$.
Hence from the definition of equivalence relation $\sim_{\eta}$ we get that
\begin{equation}
 M(P_i, \eta^i) = (\TT^{n-1} \times P_i)/\sim_{\eta} ~ \mbox{for} ~ i = 1, 2, 3.
\end{equation}


Let $U_i$ be the open subset of $\bigtriangleup^n_Q$ obtained by deleting all faces $F$ of $\bigtriangleup^n_Q$
such that the intersection $F \cap P_i$ is empty. Then $$\bigtriangleup^n_Q = U_1 \cup U_2 \cup U_3.$$ 
The space $U_i$ is diffeomorphic as manifold with corners to $ [0,1) \times P_i$. Let
$$\mathfrak{f}_i : U_i \to [0,1) \times P_i$$ be a diffeomorphism.
From the definition of $\eta$ and $\sim_{\eta}$ we get the following homeomorphisms
\begin{equation}
(\TT^{n-1} \times \mathfrak{f}_i^{-1}(\{a\} \times P_i))/\sim ~ \cong \{a\} \times M(P_i,\eta^i)
~ \mbox{for all} ~ a \in [0,1).
\end{equation}

Hence the space $\mathfrak{p}^{-1}(U_i)$ is homeomorphic to
$$(\TT^{n-1} \times \mathfrak{f}_i^{-1}([0,1) \times P_i))/\sim_{\eta} ~ \cong ~ [0,1) \times M(P_i, \eta^i).$$
Since $$W(\bigtriangleup^n_Q, \eta) = \mathfrak{p}^{-1}(U_1) \cup \mathfrak{p}^{-1}(U_2) \cup \mathfrak{p}^{-1}(U_3),$$
the space $W(\bigtriangleup^n_Q, \eta)$ is a manifold with quasitoric boundary.
The intersections $P_1\cap P_2$, $P_2\cap P_3$ and $P_1 \cap P_3$ are empty. 
Hence the boundary $$\partial{W(\bigtriangleup^n_Q, \eta)} = M(P_1,\eta^1) \sqcup M(P_2,\eta^2) \sqcup M(P_3, \eta^3).$$


\section{Orientability of $W(\bigtriangleup^n_Q, \eta)$}\label{ori2}
Fix the standard orientation on $\TT^{n-1}$. Then the boundary orientations on $P_1, P_2$ and $P_3$ induced from
the orientation of $ \bigtriangleup^n_Q $ give the orientations of $M(P_1, \eta^1), M(P_2, \eta^2)$ and
$M(P_3, \eta^3)$ respectively.

Let $W := W(\bigtriangleup^n_Q, \eta)$. The boundary $ \partial{W} $ has a collar neighborhood in $W$.
Hence by the Proposition $2.22$ of \cite{Hat}, $$H_i(W,\partial W) = \widetilde{H}_i(W/\partial W)~ \mbox{for all}~ i.$$
We show the space $W/\partial{W}$ has a $CW$-structure. Assuming $\bigtriangleup^n_Q \subset \RR^{n}$,
we choose a linear functional $$\zeta: \RR^n \to \RR$$ which distinguishes the vertices of $\bigtriangleup^n_Q$,
as in the proof of Theorem $3.1$ in \cite{DJ}. The vertices are linearly ordered according to
ascending value of $\zeta$. We make the $1$-skeleton of $\bigtriangleup^n_Q$ into a directed graph
by orienting each edge such that $\zeta$ increases along edges. For each vertex $v$ of $ \bigtriangleup^n_Q $
define its index $ ind(v) $ as the number of incident edges that point towards $ v $. Suppose
$V(\bigtriangleup^n_Q)$ is the set of vertices and $E(\bigtriangleup^n_Q)$ is the
set of edges of $\bigtriangleup^n_Q$. For each $j \in \{1, 2, \ldots, n\}$, let
$$I_j = \{(v,e_v) \in V(\bigtriangleup^n_Q) \times E(\bigtriangleup^n_Q) : ind(v)=j ~ and~
 e_v ~ is ~the~ incident ~edge~ that$$
$$points~towards~ v~ such~ that~ e_v = e \cap \bigtriangleup^n_Q ~ for ~ an ~ edge ~  e ~of ~ \bigtriangleup^n\}.$$

Suppose $(v, e_v) \in I_j$ and $F_v \subset \bigtriangleup^n_Q$ is the smallest face which contains the inward
pointing edges incident to $v$. Clearly $ind(v)=$ dim$(F_v)$. Let $U_{e_v}$ be the open subset (in relative topology)
of $F_v$ obtained by deleting all faces not containing the edge $e_v$ in $F_v$. The restriction of the equivalence
relation $\sim_{\eta}$ on $ (\TT^{n-1}\times U_{e_v})$ gives that the quotient space $(\TT^{n-1}\times U_{e_v})/\sim_{\eta}$ is
homeomorphic to the open ball $B^{2j-1} \subset \RR^{2j-1}$.

Hence the quotient space $(W/\partial W)$ has a $CW$-complex structure with odd dimensional cells and one zero
dimensional cell only. The number of $(2j-1)$-dimensional cell is $|I_j|$, the cardinality of $I_j$ for
$j=1, 2, \ldots, n$. So we get the following theorem.
\begin{theorem}
$
H_i(W, \partial W) = \left\{ \begin{array}{ll} \displaystyle \bigoplus_{|I_j|} \ZZ & \mbox{if} ~ i= 2j-1~ \mbox{and} ~ j \in \{1, \ldots, n\} \\
 \ZZ & \mbox{if}~ i = 0 \\
 0 & \mbox{otherwise}
\end{array} \right.
$
\end{theorem}
When $j=n$ the cardinality of $I_j$ is one. So $H_{2n-1}(W, \partial W) = \ZZ$. Hence  $W$ is an oriented
manifold with quasitoric boundary.
From the definition \ref{defeqiv} we get that the boundary orientation on $M(P_i, \eta^i)$ is same as the
orientation on $M(P_i, \eta^i)$ as the quasitoric manifold, for all $i = 1, 2, 3$.

\section{Identification of $M(P_1, \eta^1)$ and $M(P_2, \eta^2)$}\label{iden}
We show the quasitoric manifolds $M(P_1, \eta^1)$ and $M(P_2, \eta^2)$ are equivariantly homeomorphic up
to an automorphism of $\TT^{n-1}$. Consider the permutation $$\rho : \{0, 1, \ldots, n\} \to \{0, 1, \ldots, n\}$$
defined by
\begin{equation}
\rho(j) = \left\{ \begin{array}{ll} n-1-j & \mbox{if}~ 0 \leq j < \frac{n}{2}-1 ~ \mbox{and}~ \frac{n}{2} < j < n\\
 n & \mbox{if}~ j = \frac{n}{2}-1\\
 \frac{n}{2} & \mbox{if} ~ j= \frac{n}{2}\\
 \frac{n}{2}-1 & \mbox{if}~ j = n.
 \end{array} \right.
\end{equation}
So $\rho$ is an even or odd permutation if $n = 4l$ or $n=4l+2$ respectively.
Define a linear automorphism $ \Phi$ on $\RR^{n+1}$ by
 \begin{equation}
  \Phi(x_0, \ldots, x_j, \ldots, x_n) = (x_{\rho(0)}, \ldots, x_{\rho(j)}, \ldots, x_{\rho(n)}).
 \end{equation}
Hence $ \Phi$ is an orientation preserving or reversing diffeomorphism if $n = 4l$ or $n=4l+2$ respectively.
From equations \ref{equp1} and \ref{equp2} it is clear that $ \Phi$ maps $P_1$ diffeomorphically onto $P_2$.
We denote the restriction of $\Phi$ on the faces of $P_1$ by $\Phi$.
Also from the equations \ref{equfp1} and \ref{equfp2} we get that $\Phi$ maps the facet $\bigtriangleup_j^{n-1} \cap P_1$
of $P_1$ diffeomorphically onto the facet $\bigtriangleup_{\rho(j)}^{n-1} \cap P_2$ of $P_2$. So
\begin{equation}\label{equphi}
 \Phi(\bigtriangleup_j^{n-1} \cap P_1) = \bigtriangleup_{\rho(j)}^{n-1} \cap P_2, ~\mbox{ for all} ~ j=0,\ldots, n.
\end{equation}

Let $\alpha_1, \ldots, \alpha_{n-1}$ be the standard basis of $\ZZ^{n-1}$ over $\ZZ$.
Let $\delta^{\prime}$ be the linear automorphism of $\ZZ^{n-1}$ defined by
\begin{equation}\label{deln}
\delta^{\prime}(\alpha_{i}) = \alpha_{n-i}~  \forall ~ i= 1, \ldots, (n-1).
\end{equation}
 Hence
\begin{equation}\label{rho}
 \delta^{\prime}(\eta_{i}) = \eta_{\rho(i)} ~ \mbox{and}~ \delta^{\prime}(\eta_{\rho(i)}) = \eta_{i} ~ \mbox{for} ~ i= 0, 1, \ldots, n.
\end{equation}

Let $\delta$ be the automorphism of $\TT^{n-1}$ induced by $\delta^{\prime}$. Hence the automorphism $\delta$ is
orientation reversing if $n = 4l$ and it is orientation preserving if $4l+2$.
From the equations \ref{char1}, \ref{equphi} and \ref{rho} we get that the following commutative diagram.
$$
\begin{CD}
 \mathcal{F}(P_1) @> \Phi >> \mathcal{F}(P_2)\\
@V{\eta^1}VV @V\eta^2VV\\
\ZZ^{n-1} @> \delta >> \ZZ^{n-1}.
\end{CD}
$$

So the diffeomorphism $\delta \times \Phi : \TT^{n-1} \times P_1 \to  \TT^{n-1} \times P_2$
induces a $\delta$-equivariant orientation reversing homeomorphism $$g_n : M(P_1, \eta^1) \to M(P_2, \eta^2).$$
From the definition \ref{char1} of the characteristic function $ \eta^{3}$ we get that the quasitoric manifold
$M(P_3, \eta^3)$ is equivariantly homeomorphic to $\CP^{n-1} $ if  $n=4l+2$ and $\overline{\CP}^{n-1} $ if $n = 4l$.
Define an equivalence relation $\sim_n$ on $ W(\bigtriangleup^n_Q, \eta) $ by
\begin{equation}
 x \sim_n y ~\mbox{if and only if} ~ x \in M(P_1, \eta^1) ~ \mbox{and} ~ y = g_n(x).
\end{equation}
So the quotient space $ W(\bigtriangleup^n_Q, \eta)/\sim_n $ is an oriented manifold with boundary. The boundary
of this manifold is $\CP^{n-1}$ if $n =4l+2$ and the boundary is $\overline{\CP}^{n-1}$ if $n=4l$.
Hence we have proved the following theorem.

\begin{theorem}\label{gr}
The complex projective space $\CP^{2k+1}$ is the boundary of an oriented manifold, for each $k \geq 0$.
\end{theorem}

\begin{example}
We adhere to definition and notations given in the example \ref{exam1}.
The faces  $A_0A_1$ and $A_3A_4$ are the intersection of facets $\{\bigtriangleup_{4}^3, \bigtriangleup_{3}^3, \bigtriangleup_2^3\}$
and $\{\bigtriangleup_{2}^3, \bigtriangleup_{1}^3, \bigtriangleup_0^3\}$ respectively of $\bigtriangleup^4$.

Here the simple polytopes $P_1, P_2$ are prism and $P_3$ is a $3$-simplex, see the figure \ref{egch08}.
The restriction of the function $\eta$ (namely $\eta^{1}$, $\eta^{2}$ and $\eta^{3}$) on the facets of
$P_1, P_2$ and $P_3$ are given in following figure \ref{eg5}.
\begin{figure}[ht]
\centerline{
\scalebox{0.60}{
\input{egch08.pstex_t}
 }
 }
\caption {The simple polytope $ \bigtriangleup^4_Q$ with the facets $P_1$, $P_2$ and $ P_3$.}
\label{egch08}
\end{figure}

\begin{figure}[ht]
\centerline{
\scalebox{0.60}{
\input{eg5.pstex_t}
 }
 }
\caption {The characteristic functions $\eta^{1}$, $\eta^{2}$ and $\eta^{3}$ of $P_1, P_2$ and $P_3$ respectively.}
\label{eg5}
\end{figure}
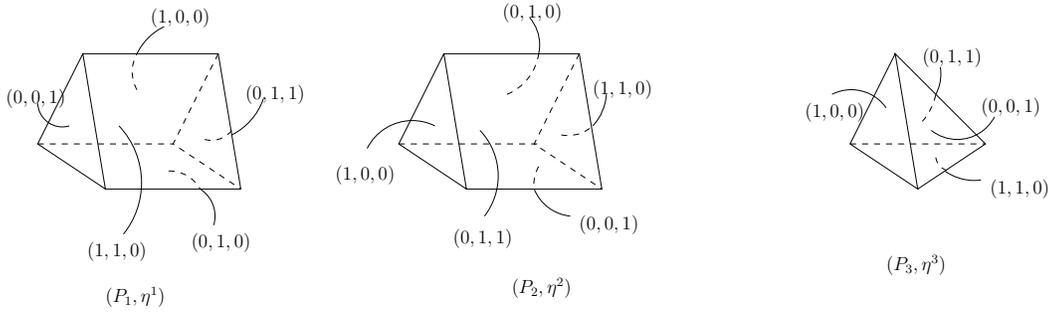
Let $\delta^{\prime}$ be the automorphism of $\ZZ^{3}$ defined by
$$\delta^{\prime}(1,0,0) = (0,0,1), \delta^{\prime}(0,1,0) = (0,1,0)~ \mbox{and}~ \delta^{\prime}(0,0,1)=(1,0,0).$$
Clearly the characteristic pairs $(P_1, \eta^{1})$ and $ (P_2, \eta^{2})$ give two $\delta$-equivariantly homeomorphic
quasitoric manifolds, namely $M(P_1, \eta^1)$ and $M(P_2, \eta^2)$ respectively.
The combinatorial pair $(P_3, \eta^{3})$ gives the quasitoric manifold $\overline{\CP}^3$ over $P_3$.
So the boundary of $W(\bigtriangleup^4_Q, \eta)$ is $M(P_1, \eta^1) \sqcup M(P_2, \eta^2) \sqcup \overline{\CP}^3$.
Hence after identifying  $M(P_1, \eta^1)$ and $M(P_2, \eta^2)$ via an orientation reversing homeomorphism,
we get an oriented manifold with boundary and the boundary is $\overline{\CP}^3$.
\end{example}

{\bf Acknowledgement.} I would like to thank my thesis supervisor, Professor Mainak Poddar, and 
Professor Goutam Mukherjee for helpful suggestions and stimulating discussions. I thank
Indian Statistical Institute for supporting my research fellowship. 


\renewcommand{\refname}{References}

\vspace{1cm}

\vfill

\end{document}

%% file: egch060.pstex_t
\begin{picture}(0,0)%
\includegraphics{egch060.pstex}%
\end{picture}%
\setlength{\unitlength}{4144sp}%
\begingroup\makeatletter\ifx\SetFigFontNFSS\undefined%
\gdef\SetFigFontNFSS#1#2#3#4#5{%
  \reset@font\fontsize{#1}{#2pt}%
  \fontfamily{#3}\fontseries{#4}\fontshape{#5}%
  \selectfont}%
\fi\endgroup%
\begin{picture}(2550,2928)(1381,-3100)
\put(1396,-1051){\makebox(0,0)[lb]{\smash{{\SetFigFontNFSS{12}{14.4}{\rmdefault}{\mddefault}{\updefault}{\color[rgb]{0,0,0}$A_0$}%
}}}}
\put(1621,-3031){\makebox(0,0)[lb]{\smash{{\SetFigFontNFSS{12}{14.4}{\rmdefault}{\mddefault}{\updefault}{\color[rgb]{0,0,0}$A_{1}$}%
}}}}
\put(3511,-3031){\makebox(0,0)[lb]{\smash{{\SetFigFontNFSS{12}{14.4}{\rmdefault}{\mddefault}{\updefault}{\color[rgb]{0,0,0}$A_{2}$}%
}}}}
\put(3916,-1861){\makebox(0,0)[lb]{\smash{{\SetFigFontNFSS{12}{14.4}{\rmdefault}{\mddefault}{\updefault}{\color[rgb]{0,0,0}$A_{3}$}%
}}}}
\put(2611,-331){\makebox(0,0)[lb]{\smash{{\SetFigFontNFSS{12}{14.4}{\rmdefault}{\mddefault}{\updefault}{\color[rgb]{0,0,0}$A_4$}%
}}}}
\end{picture}%

%% file: egch06.pstex_t
\begin{picture}(0,0)%
\includegraphics{egch06.pstex}%
\end{picture}%
\setlength{\unitlength}{4144sp}%
\begingroup\makeatletter\ifx\SetFigFontNFSS\undefined%
\gdef\SetFigFontNFSS#1#2#3#4#5{%
  \reset@font\fontsize{#1}{#2pt}%
  \fontfamily{#3}\fontseries{#4}\fontshape{#5}%
  \selectfont}%
\fi\endgroup%
\begin{picture}(2685,3411)(1246,-3628)
\put(1441,-961){\makebox(0,0)[lb]{\smash{{\SetFigFontNFSS{12}{14.4}{\rmdefault}{\mddefault}{\updefault}{\color[rgb]{0,0,0}$A_{1}$}%
}}}}
\put(2566,-376){\makebox(0,0)[lb]{\smash{{\SetFigFontNFSS{12}{14.4}{\rmdefault}{\mddefault}{\updefault}{\color[rgb]{0,0,0}$A_{0}$}%
}}}}
\put(1621,-3031){\makebox(0,0)[lb]{\smash{{\SetFigFontNFSS{12}{14.4}{\rmdefault}{\mddefault}{\updefault}{\color[rgb]{0,0,0}$A_{2}$}%
}}}}
\put(3511,-3031){\makebox(0,0)[lb]{\smash{{\SetFigFontNFSS{12}{14.4}{\rmdefault}{\mddefault}{\updefault}{\color[rgb]{0,0,0}$A_{3}$}%
}}}}
\put(3916,-1861){\makebox(0,0)[lb]{\smash{{\SetFigFontNFSS{12}{14.4}{\rmdefault}{\mddefault}{\updefault}{\color[rgb]{0,0,0}$A_{4}$}%
}}}}
\put(1261,-2401){\makebox(0,0)[lb]{\smash{{\SetFigFontNFSS{12}{14.4}{\rmdefault}{\mddefault}{\updefault}{\color[rgb]{0,0,0}$H_3$}%
}}}}
\put(3151,-3526){\makebox(0,0)[lb]{\smash{{\SetFigFontNFSS{12}{14.4}{\rmdefault}{\mddefault}{\updefault}{\color[rgb]{0,0,0}$H_2$}%
}}}}
\put(3331,-736){\makebox(0,0)[lb]{\smash{{\SetFigFontNFSS{12}{14.4}{\rmdefault}{\mddefault}{\updefault}{\color[rgb]{0,0,0}$H_1$}%
}}}}
\end{picture}%

%% file: egch08.pstex_t
\begin{picture}(0,0)%
\includegraphics{egch08.pstex}%
\end{picture}%
\setlength{\unitlength}{4144sp}%
\begingroup\makeatletter\ifx\SetFigFontNFSS\undefined%
\gdef\SetFigFontNFSS#1#2#3#4#5{%
  \reset@font\fontsize{#1}{#2pt}%
  \fontfamily{#3}\fontseries{#4}\fontshape{#5}%
  \selectfont}%
\fi\endgroup%
\begin{picture}(3897,3816)(1966,-5665)
\put(1981,-4336){\makebox(0,0)[lb]{\smash{{\SetFigFontNFSS{12}{14.4}{\rmdefault}{\mddefault}{\updefault}{\color[rgb]{0,0,0}$P_1$}%
}}}}
\put(4951,-5596){\makebox(0,0)[lb]{\smash{{\SetFigFontNFSS{12}{14.4}{\rmdefault}{\mddefault}{\updefault}{\color[rgb]{0,0,0}$P_3$}%
}}}}
\put(5446,-2671){\makebox(0,0)[lb]{\smash{{\SetFigFontNFSS{12}{14.4}{\rmdefault}{\mddefault}{\updefault}{\color[rgb]{0,0,0}$P_2$}%
}}}}
\end{picture}%

%% file: eg5.pstex_t
\begin{picture}(0,0)%
\includegraphics{eg5.pstex}%
\end{picture}%
\setlength{\unitlength}{4144sp}%
\begingroup\makeatletter\ifx\SetFigFontNFSS\undefined%
\gdef\SetFigFontNFSS#1#2#3#4#5{%
  \reset@font\fontsize{#1}{#2pt}%
  \fontfamily{#3}\fontseries{#4}\fontshape{#5}%
  \selectfont}%
\fi\endgroup%
\begin{picture}(9878,3063)(571,-3505)
\put(1396,-2986){\makebox(0,0)[lb]{\smash{{\SetFigFontNFSS{12}{14.4}{\rmdefault}{\mddefault}{\updefault}{\color[rgb]{0,0,0}$(1, 1, 0)$}%
}}}}
\put(2431,-2896){\makebox(0,0)[lb]{\smash{{\SetFigFontNFSS{12}{14.4}{\rmdefault}{\mddefault}{\updefault}{\color[rgb]{0,0,0}$(0,1,0)$}%
}}}}
\put(2026,-646){\makebox(0,0)[lb]{\smash{{\SetFigFontNFSS{12}{14.4}{\rmdefault}{\mddefault}{\updefault}{\color[rgb]{0,0,0}$(1,0,0)$}%
}}}}
\put(586,-1456){\makebox(0,0)[lb]{\smash{{\SetFigFontNFSS{12}{14.4}{\rmdefault}{\mddefault}{\updefault}{\color[rgb]{0,0,0}$(0,0,1)$}%
}}}}
\put(6301,-2716){\makebox(0,0)[lb]{\smash{{\SetFigFontNFSS{12}{14.4}{\rmdefault}{\mddefault}{\updefault}{\color[rgb]{0,0,0}$(0,0,1)$}%
}}}}
\put(3871,-2221){\makebox(0,0)[lb]{\smash{{\SetFigFontNFSS{12}{14.4}{\rmdefault}{\mddefault}{\updefault}{\color[rgb]{0,0,0}$(1,0,0)$}%
}}}}
\put(8551,-1591){\makebox(0,0)[lb]{\smash{{\SetFigFontNFSS{12}{14.4}{\rmdefault}{\mddefault}{\updefault}{\color[rgb]{0,0,0}$(1,0,0)$}%
}}}}
\put(10306,-1546){\makebox(0,0)[lb]{\smash{{\SetFigFontNFSS{12}{14.4}{\rmdefault}{\mddefault}{\updefault}{\color[rgb]{0,0,0}$(0,0,1)$}%
}}}}
\put(2971,-1411){\makebox(0,0)[lb]{\smash{{\SetFigFontNFSS{12}{14.4}{\rmdefault}{\mddefault}{\updefault}{\color[rgb]{0,0,0}$(0,1,1)$}%
}}}}
\put(6436,-1366){\makebox(0,0)[lb]{\smash{{\SetFigFontNFSS{12}{14.4}{\rmdefault}{\mddefault}{\updefault}{\color[rgb]{0,0,0}$(1,1,0)$}%
}}}}
\put(5536,-601){\makebox(0,0)[lb]{\smash{{\SetFigFontNFSS{12}{14.4}{\rmdefault}{\mddefault}{\updefault}{\color[rgb]{0,0,0}$(0,1,0)$}%
}}}}
\put(5041,-2851){\makebox(0,0)[lb]{\smash{{\SetFigFontNFSS{12}{14.4}{\rmdefault}{\mddefault}{\updefault}{\color[rgb]{0,0,0}$(0,1,1)$}%
}}}}
\put(9721,-1051){\makebox(0,0)[lb]{\smash{{\SetFigFontNFSS{12}{14.4}{\rmdefault}{\mddefault}{\updefault}{\color[rgb]{0,0,0}$(0,1,1)$}%
}}}}
\put(10396,-2356){\makebox(0,0)[lb]{\smash{{\SetFigFontNFSS{12}{14.4}{\rmdefault}{\mddefault}{\updefault}{\color[rgb]{0,0,0}$(1,1,0)$}%
}}}}
\put(1576,-3436){\makebox(0,0)[lb]{\smash{{\SetFigFontNFSS{12}{14.4}{\rmdefault}{\mddefault}{\updefault}{\color[rgb]{0,0,0}$( P_1, \eta^{1})$}%
}}}}
\put(5626,-3346){\makebox(0,0)[lb]{\smash{{\SetFigFontNFSS{12}{14.4}{\rmdefault}{\mddefault}{\updefault}{\color[rgb]{0,0,0}$( P_2, \eta^{2})$}%
}}}}
\put(9361,-3121){\makebox(0,0)[lb]{\smash{{\SetFigFontNFSS{12}{14.4}{\rmdefault}{\mddefault}{\updefault}{\color[rgb]{0,0,0}$( P_3, \eta^{3})$}%
}}}}
\end{picture}%